\documentclass[12pt]{amsart}
\usepackage{fullpage}
\usepackage{amssymb}
\usepackage{graphicx}

\newtheorem{theorem}{Theorem!!!}[section]
\newtheorem{thm}[theorem]{Theorem}

\newtheorem{prop}[theorem]{Proposition}
\newtheorem{lem}[theorem]{Lemma}

\newtheorem{conj}[theorem]{Conjecture}

\theoremstyle{definition}
\newtheorem{defn}{Definition}[section]

\newcommand{\thmref}[1]{Theorem~\ref{#1}}

\numberwithin{equation}{section}

\newcommand{\propref}[1]{Proposition~\ref{#1}}
\newcommand{\secref}[1]{\S\ref{#1}}
\newcommand{\lemref}[1]{Lemma~\ref{#1}}

\newcommand{\beq}{\begin{equation}}
\newcommand{\eeq}{\end{equation}}
\newcommand{\eref}[1]{(\ref{#1})}

\newcommand{\comm}[1]{}

\newcommand{\BBB}[1]{{\mathbb #1}}
\newcommand{\ra}{\rightarrow}

\newcommand{\de}{\delta}
\newcommand{\dist}{\operatorname{dist}}
\newcommand{\al}{\alpha}
\newcommand{\be}{\beta}
\newcommand{\eps}{\epsilon}

\newcommand{\ga}{\gamma}
\newcommand{\Om}{\Omega}
\newcommand{\om}{\omega}

\newcommand{\ve}{\varepsilon}

\newcommand{\cB}{{\mathcal B}}

\newcommand{\cC}{{\mathcal C}}

\newcommand{\NN}{{\BBB N}}
\newcommand{\TT}{{\BBB T}}
\newcommand{\ZZ}{{\BBB Z}}
\newcommand{\QQ}{{\BBB Q}}

\newcommand{\RR}{{\BBB R}}

\newcommand{\CC}{{\BBB C}}
\newcommand{\DD}{{\BBB D}}

\begin{document}

\title[Constructing Non-computable Julia sets]{Constructing Non-Computable Julia Sets}

\author{M. Braverman}
\author{M. Yampolsky}
\date{\today}
\begin{abstract}
We completely characterize the conformal radii of Siegel disks
in the family $$P_\theta(z)=e^{2\pi i\theta}z+z^2,$$
corresponding to {\bf computable} parameters $\theta$. As a consequence, we constructively produce
quadratic polynomials with {\bf non-computable} Julia sets.

\end{abstract}
\maketitle

The purpose of this note is to completely characterize the conformal radii of Siegel disks
in the family $$P_\theta(z)=e^{2\pi i\theta}z+z^2,$$
corresponding to  computable parameters $\theta$. 
As one consequence, we derive the following statement:

\medskip
\noindent
{\bf Theorem.} {\it There exist computable complex parameters $c$, such that the Julia
set $J_c$ of the quadratic polynomial $f_c(z)=z^2+c$ is non-computable.}

\medskip
\noindent
In the end of the note we discuss the computational complexity of constructing such
values of $c$, and show that it is poly-time, assuming a certain 
widely believed conjecture holds.

\section{Computability of subsets of $\RR^n$}
We recall the relevant definitions of constructive analysis very briefly. The reader is 
referred to our earlier work \cite{BY} for a more detailed exposition.

Denote by $\DD$ the set of the {\it dyadic rationals}, that is, rationals of 
the form $\frac{p}{2^m}$. 
The classical definition of a computable real number may be formulated as follows:
\begin{defn}
A number $x\in\RR$ is computable, if there exists a Turing Machine (further abbreviated
as TM) $M(n)$ with a positive integer input $n$, such that for each $n$, 
this TM terminates and outputs a number $d_n\in\DD$ with the property $|x-d_n|<2^{-n}$.
\end{defn}

\noindent
The definition is generalized to points in $\RR^k$ in an obvious way.

Sometimes, it is desirable to use a real number as a parameter of a
computation, without regard to its computable properties. This is 
achieved with the use of oracles.

\begin{defn}
We say that $\phi: \NN \rightarrow \DD$ is an {\it oracle}
for a real number $x$, if $| x - \phi(n)|<2^{-n}$ for all $n \in \NN$. 
We say that a TM $M^{\phi}$ is 
an {\it oracle machine}, if at every step of 
the computation M is allowed to query the value $\phi(n)$ for any $n$.
\end{defn}

\noindent
Now, for instance, we may define computable functions of a real variable
(compare with the discussion in \cite{Brv05}).

\begin{defn}
\label{funcomp}
We say that a function $f:[a,b] \rightarrow [c,d]$ is computable, if 
there exists an oracle TM $M^{\phi} (m)$ such that if $\phi$ is an oracle 
for $x \in [a,b]$, then on input $m$, $M^{\phi}$ outputs a $y \in \DD$ 
such that $| y - f(x)|<2^{-m}$.
\end{defn}

\noindent
Recall that the Hausdorff metric is a metric on 
compact subsets of $\RR^n$ defined by 
\begin{equation}
\label{hausdorff metric}
d_H ( X, Y) =  \inf \{\epsilon > 0 \;|\; X \subset U_{\epsilon} 
(Y)~~\mbox{and}~~  Y \subset U_{\epsilon}(X)\},
\end{equation}

\noindent
where $U_\eps(S)$ is defined as the union of the set of $\eps$-balls with centers in $S$.

\noindent
 We introduce
a class $\cC$ of sets which is dense in metric $d_H$ among the compact 
sets and which has a natural correspondence to binary strings. 
Namely $\cC$ is the set of finite unions of dyadic balls:
$$
\cC= \left\{ \bigcup_{i=1}^n \overline{B(d_i, r_i)}~|~~\mbox{where}~~d_i, 
r_i 
\in \DD \right\}.
$$

\noindent
We now define the notion of computability of subsets of $\RR^n$ (see \cite{Wei}).

\begin{defn}
\label{setcomp}
We say that a compact set $K\Subset \RR^k$ is computable, if there exists a
TM $M(n)$, $n\in \NN$ which outputs a set $C_n\in\cC$ such that
$\dist_H(C_n,K)<2^{-n}.$

\end{defn}

\noindent
It is not difficult to see that:

\begin{prop}
A set $K\Subset \RR^k$ is computable if and only if the {\em distance function} $d_K (x) = \inf \{ |x-y|~~|~~y\in K \}$
is a computable function of a real variable.
\end{prop}

In this note we are concerned with computability of quadratic Julia sets
$J_c  = J(z^2+c)$.
Define the function $J: \CC \rightarrow K^{*}$ ($K^{*}$ is 
the set of all compact subsets of $\CC$) by $J(c)=J(f_c)$. 
In a complete analogy to Definition \ref{funcomp} we can define

\begin{defn}
\label{funcomp2}
We say that a function $\kappa:S \rightarrow K^{*}$ for some bounded set $S$ is 
computable, if there exits  an oracle TM $M^\phi(m)$ with $\phi$  representing
$x\in S$ such that 
 on input $m$, $M^{\phi}$ outputs a $C \in \cC$ 
such that $d_H (C, \kappa(x))<2^{-m}$.
\end{defn}

\noindent
In the case of Julia sets:

\begin{defn}
\label{Jcomp}
We say that $J_c$ is computable if the function $J: d \mapsto J_d$ is 
computable on the set $\{ c \}$.
\end{defn}

\noindent
We have established the existence of non-computable quadratic Julia sets 
in \cite{BY}, based on the analysis of Julia sets with Siegel disks. The
next chapter summarizes the relevant tools of Complex Dynamics.

\section{Conformal radii of Siegel disks and Yoccoz's Brjuno function}

Let $R:\hat\CC\to\hat\CC$ be a rational map of the Riemann sphere.
For a periodic point $z_0=R^p(z_0)$
of period $p$ its {\it multiplier} is the quantity $\lambda=\lambda(z_0)=DR^p(z_0)$.
We may speak of the multiplier of a periodic cycle, as it is the same for all points
in the cycle by the Chain Rule. In the case when $|\lambda|\neq 1$, the dynamics
in a sufficiently small neighborhood of the cycle is governed by the Intermediate 
Value Theorem: when $0<|\lambda|<1$, the cycle is {\it attracting} ({\it super-attracting}
if $\lambda=0$), if $|\lambda|>1$ it is {\it repelling}.
Both in the attracting and repelling cases, the dynamics can be locally linearized:
\begin{equation}
\label{linearization-equation}
\psi(R^p(z))=\lambda\cdot\psi(z)
\end{equation}
where $\psi$ is a conformal mapping of a small neighborhood of $z_0$ to a disk around $0$.
By a classical result of Fatou, a rational mapping has at most finitely many non-repelling
periodic orbits.

In the case when $\lambda=e^{2\pi i\theta}$, $\theta\in\RR$, 
 the simplest to study is the {\it parabolic case} when $\theta=n/m\in\QQ$, so $\lambda$ 
is a root of unity. In this case $R^p$ is not locally linearizable; it is not hard to see that $z_0\in J(R)$.
 In the complementary situation, two non-vacuous possibilities  are considered:
{\it Cremer case}, when $R^p$ is not linearizable, and {\it Siegel case}, when it is.
In the latter case, the linearizing map $\psi$ from (\ref{linearization-equation}) conjugates
the dynamics of $R^p$ on a neighborhood $U(z_0)$ to the irrational rotation by angle $\theta$
(the {\it rotation angle})
on a disk around the origin. The maximal such neighborhood of $z_0$ is called a {\it Siegel disk}.

Let us discuss in more detail the occurrence of Siegel disks in the quadratic family.
For a number $\theta\in [0,1)$ denote $[r_0,r_1,\ldots,r_n,\ldots]$, $r_i\in\NN\cup\{\infty\}$ its possibly finite 
continued fraction expansion:
\begin{equation}
\label{cfrac}
[r_0,r_1,\ldots,r_n,\ldots]\equiv\cfrac{1}{r_0+\cfrac{1}{r_1+\cfrac{1}{\cdots+\cfrac{1}{r_n+\cdots}}}}
\end{equation}
Such an expansion is defined uniquely if and only if $\theta\notin\QQ$. In this case, the {\it rational 
convergents } $p_n/q_n=[r_0,\ldots,r_{n-1}]$ are the closest rational approximants of $\theta$ among the
numbers with denominators not exceeding $q_n$. In fact, setting $\lambda=e^{2\pi i\theta}$, we have
$$|\lambda^h-1|>|\lambda^{q_n}-1|\text{ for all }0<h<q_{n+1},\; h\neq q_n.$$
The difference $|\lambda^{q_n}-1|$ lies between $2/q_{n+1}$ and $2\pi/q_{n+1}$,
therefore the rate of growth of the denominators $q_n$ describes how well 
$\theta$ may be approximated with rationals.

We recall a theorem due to Brjuno (1972):
\begin{thm}[\cite{Bru}]
Let $R$ be an analytic map with a periodic point $z_0\in\hat\CC$. Suppose that the 
multiplier of $z_0$ is $\lambda=e^{2\pi i\theta}$, and
\begin{equation}
\label{brjuno}
B(\theta)=\displaystyle\sum_n\frac{\log(q_{n+1})}{q_n}<\infty.
\end{equation}
Then $z_0$ is a Siegel point.
\end{thm}

\noindent
Note that a quadratic polynomial with a fixed Siegel disk with rotation angle $\theta$ after an affine
change of coordinates can be written as 
\begin{equation}
\label{form-1}
P_\theta(z)=z^2+e^{2\pi i \theta}z.
\end{equation}
\noindent
In 1987 Yoccoz \cite{Yoc} proved the following converse to Brjuno's Theorem:

\begin{thm}[\cite{Yoc}]
Suppose that for $\theta\in[0,1)$ the polynomial $P_\theta$ has a Siegel point at the origin.
Then $B(\theta)<\infty$.
\end{thm}

\noindent
The numbers satisfying (\ref{brjuno}) are called Brjuno numbers; the set of all Brjuno numbers will be denoted $\cB$. It is a full measure set which contains all Diophantine rotation numbers.
In particular, the rotation numbers $[r_0,r_1,\ldots]$ of {\it bounded type},
that is with $\sup r_i<\infty$ are in $\cB$.
The sum of the series (\ref{brjuno}) is called the Brjuno function. 
For us a different characterization of $\cB$ will be more useful. Inductively define $\theta_1=\theta$
and $\theta_{n+1}=\{1/\theta_n\}$. In this way, 
$$\theta_n=[r_{n-1},r_n,r_{n+1},\ldots].$$
We define the {\it Yoccoz's Brjuno function} as
$$\Phi(\theta)=\displaystyle\sum_{n=1}^{\infty}\theta_1\theta_2\cdots\theta_{n-1}\log\frac{1}{\theta_n}.$$
One can verify that $$B(\theta)<\infty\Leftrightarrow \Phi(\theta)<\infty.$$
The value of the function $\Phi$ is related to the size of the Siegel disk in the following way.

\begin{defn}
Let $(U,u)$ be a simply-connected subdomain of $\CC$ with a marked interior point.
Consider the unique conformal isomorphism $\phi:\DD\mapsto U$ with
$\phi(0)=u$, and $\phi'(0)>0$. The {\it conformal radius of $(U,u)$} is the value of the
derivative $r(U,u)=\phi'(0)$.

Let $P(\theta)$ be a quadratic polynomial with a Siegel disk $\Delta_\theta\ni 0$. 
The {\it conformal radius of the Siegel disk $\Delta_\theta$} is
$r(\theta)=r(\Delta_\theta,0)$.
For all other $\theta\in[0,\infty)$ we set $r(\theta)=0$, and $\Delta_\theta=\{0\}$. 
\end{defn} 

\noindent
By the Koebe 1/4 Theorem of classical complex analysis (see e.g. \cite{Ahlfors}), the  radius of 
the largest Euclidean disk around $u$ which can be inscribed in $U$ is at least $r(U,u)/4$.

\noindent
We note that one has the following direct consequence of the Carath{\'e}odory Kernel Theorem
(see e.g. \cite{Pom}):

\begin{prop}
\label{radius-continuous}
The conformal radius of a quadratic Siegel disk varies continuously with respect to the Hausdorff 
distance on Julia sets.
\end{prop}

\noindent
Yoccoz \cite{Yoc} has shown that the sum 
$$\Phi(\theta)+\log r(\theta)$$
is bounded below independently of $\theta\in\cB$. Recently, Buff and Ch{\'e}ritat have greatly improved this result
by showing that:

\begin{thm}[\cite{BC}]
\label{phi-cont}
The function $\theta\mapsto \Phi(\theta)+\log r(\theta)$ extends to $\RR$ as a 1-periodic continuous
function.
\end{thm}

\medskip
\noindent
In \cite{BBY1} we obtain the following result on 
computability of quadratic Siegel disks:

\begin{thm}
\label{thm BBY}
The following statements are equivalent:
\begin{itemize}
\item[(I)] the Julia set $J(P_\theta)$ is computable by a TM with an oracle access to $\theta$;
\item[(II)] the conformal radius $r(\theta)$ is computable by a TM with an oracle access to $\theta$;
\item [(III)] the inner radius $\inf_{z\in\partial \Delta_\theta}|z|$ is computable by a TM with an oracle access to $\theta$.
\end{itemize}
\end{thm}

\noindent
We note that when $\theta$ is not a Brjuno number, the quantities in (II) and (III) are each
equal to zero, and the claim is simply that $J(P_\theta)$ is computable in this case.

\medskip

\noindent
We will make use of the following Lemma which bounds
the variation of the conformal radius under a perturbation of the domain.
It  is a direct consequence of the Koebe Theorem (see e.g. \cite{RZ} for a
proof).

\begin{lem}
\label{rad-modulus1}
Let $U$ be a simply-connected subdomain of $\CC$ containing the point $0$ in the
interior. Let $V \subset U$ be a subdomain of $U$. 
Assume that $\partial V\subset B_\eps(\partial U)$.
Then 
$$0<r(U,0)-r(V,0)\leq 4 \sqrt{r(U, 0)}\sqrt{\eps}.$$
\end{lem}

\section{Computing noble Siegel discs}
We will make use of a computability result for {\it noble} Siegel disks.
The term ``noble'' is applied in the literature to rotation numbers
of the form $[a_0,a_1,\ldots,a_k,1,1,1,\ldots]$. The noblest of 
all is the golden mean $\gamma_*=[1,1,1,\ldots]$.

\begin{lem}
\label{bounded-type}
There is a Turing Machine $M$, which given a finite sequence
of numbers $[a_0, a_1, \ldots, a_k]$ computes the conformal radius
$r_\gamma$ for 
the noble number $\gamma=[a_0,  \ldots, a_k, 1,  \ldots]$. 
\end{lem}

\medskip
\noindent

\comm{
The idea is to approximate the boundary of $\Delta_{\gamma}$
with the iterates of the critical point $c_{\gamma}=-e^{2\pi i\gamma}/2$.
It is known that in this case the critical point itself is contained in the boundary.
The renormalization theory for golden-mean Siegel disks (constructed in \cite{McM}) implies
that the boundary $\Delta_{\gamma_*}$ is self-similar up to an 
exponentially small error. In particular, there exist  constants
$C>0$ and $\lambda>1$ such 
that
\begin{equation*}
d_H(\{P_{\gamma_*}^i(c_{\gamma_*}),i=0,\ldots, q_n\},\partial \Delta_{\gamma_*})<C\lambda^{-n}
\end{equation*}

Below we derive a similar estimate for all noble Siegel disks with 
constructive constants $C$ and $\lambda$.
For this, we do not need to invoke the whole power of renormalization theory. Rather,
we will use a theorem of Douady, Ghys, Herman, and Shishikura \cite{Do1} which specifically
applies to quadratic noble Siegel disks.

}

Noble (or more generally, bounded type) Sigel quadratic Julia sets may be 
constructed by means of quasiconformal surgery (cf. \cite{Do1}) on a Blaschke product
$$f_\gamma(z)=e^{2\pi i\tau(\gamma)}z^2\frac{z-3}{1-3z}.$$
This map homeomorphically maps the unit circle $\TT$ onto itself with
a single (cubic) critical point at $1$. The angle
$\tau(\gamma)$ can be uniquely selected in such a way that the rotation number
of the restriction $\rho(f_\gamma|_\TT)=\gamma$.

For each $n$, the points $$\{1, f_\gamma(1),f^2_\gamma(1),\ldots,f^{q_{n+1}-1}_\gamma(1)\}$$
form the {\it $n$-th dynamical partition} of the unit circle. We have (cf. Theorem 3.1 of \cite{dFdM})
the following:

\begin{thm}[{\bf Universal real {\it a priori} bound}]
\label{real bound}
There exists an explicit constant $B>1$ independent of $\gamma$ and $n$
such that the following holds. 
Let $\gamma\in\RR\setminus\QQ$ and $n\in\NN$. Then
any two adjacent intervals $I$ and $J$
of the $n$-th dynamical partition of $f_\gamma$ are $B$-commensurable:
$$B^{-1}|I|\leq |J|\leq B|I|.$$
\end{thm}

\begin{prop}[\cite{He}]
\label{qs-conjugacy}
For each noble $\gamma=[a_0,  \ldots, a_k, 1,  \ldots]$ the Blaschke product $f_\gamma$
is $K_1$-quasisymmetrically conjugate to the rotation $R_\gamma:x\mapsto x+\gamma \mod \ZZ$.
The quasisymmetric constant may be taken as $K_1={(2\max a_i)}^{10B^2}$.
\end{prop}

\noindent
Let us now consider the mapping 
$\Psi$ which identifies the critical orbits of $f_\gamma$ and $P_\gamma$ by
$$\Psi:f^i_\gamma(1)\mapsto P^i_\gamma(c_\gamma).$$

\noindent
We have the following (see, for example, Theorem 3.10 of \cite{YZ}):

\begin{thm}[{\bf Douady, Ghys, Herman, Shishikura}]
\label{surgery}
The mapping $\Psi$ extends to a $K$-quasiconformal homeomorphism of the plane
$\CC$ which maps the unit disk $\DD$ onto the Siegel disk $\Delta_\gamma$.
The constant $K$ may be taken as the quasiconformal dilatation of 
any global quasiconformal extension of the $K_1$-qs conjugacy of \propref{qs-conjugacy}.
In particular, $K\leq 2K_1$.
\end{thm}

Elementary combinatorics implies that each interval of the $n$-th dynamical
partition contains at least two intervals of the $(n+2)$-nd dynamical partition.
This in conjunction with \thmref{real bound} implies that the size of an interval
of the $(n+2)$-nd dynamical partition of $f_\gamma$ is at most $\tau^n$ where
$$\tau=\sqrt\frac{B}{B+1}.$$

We now  complete the proof of \lemref{bounded-type}. Denote $W_n$ the connected component
containing $0$ of the domain obtained by removing from the plane a closed disk of
radius $2K\tau^n$ around each point of 
$$\Omega_n=\{P_\gamma^i(c_\gamma),\;i=0,\ldots,q_{n+2}\}.$$
By Theorem \ref{surgery}, 
$$ \dist_H ( \Omega_n, \partial \Delta_\gamma) < K \tau^n,  $$
and we have $$W_n\subset \Delta_\gamma\text{ and }\dist_H(\partial \Delta_\gamma,\partial W_n)\leq 
\eps_n=2K\tau^n.$$
Any constructive algorithm for producing the Riemann mapping of a planar region 
(e.g. that of \cite{BB}) can be used to estimate the conformal radius $r(W_n,0)$ with
precision $\eps_n$. Denote this estimate $r_n$.

Elementary estimates imply that the Julia set $J(P_\gamma)\subset \overline{B(0,2)}$. By
Schwarz Lemma this implies 
$r(\Delta_\gamma,0)<2$.
By \lemref{rad-modulus1} we have 
$$|r(\Delta_\gamma,0)-r_n|
\le |r(\Delta_\gamma,0)- r(W_n,0)| + | r(W_n,0) - r_n | 
<4\sqrt{\eps_n}+\eps_n\underset{n\to\infty}{\longrightarrow}0,$$
and the proof is complete.

\section{Main result}
\begin{defn}(cf. \cite{Wei})
A real number $r$ is {\it right-computable} (or {\it right r.e.})
 if there exists a TM $M(n)$ 
which outputs a non-increasing sequence of dyadic numbers $r_n\geq r$ 
such that $r_n\searrow r$.
\end{defn}
\noindent
A right computable number does not have to be computable, as the definition
does not require the existence of a constructive bound on the difference
$r_n-r$. We will give a specific example of a non-computable $r$ in \secref{sec:constructive}.

The following result gives a precise characterization of the conformal 
radii of quadratic Siegel Julia sets $J_{z^2+c}$ with a computable $c$. 

\begin{thm}
\label{thm:char} 
Let $r \in [0,0.1]$ be a real number. Then $r=r(c)$ is the conformal 
radius of a Siegel disc of the Julia set $J_{z^2+c}$ for some computable 
number $c$ if and only if $r$ is right-computable.
\end{thm}

\begin{proof}
We prove the ``only if" direction here. The remainder of 
the section is dedicated to proving the ``if" direction. 

We assume that $c$ is computable, and show that $r(c)$ is right-computable. 
The case $r(c)=0$ is obvious, we assume that $r(c)>0$ for the remainder
of the proof. 
Recall that periodic orbits are dense in the Julia set $J_c$. Let $H_n$ 
be the union of the first $n$ repelling periodic orbits. Let $\al$ be 
the center of the Siegel disc we are considering. Note that $\al$ is a computable 
real number.

We can algorithmically find a strictly increasing sequence $\{n_l\}\subset\NN$
such that $\CC\setminus U_l$ has a simply-connected component $W_l$ containing $\al$,
where 
$$
B_{2^{-(l+1)}}(H_{n_l})\subset U_l\subset B_{2^{-l}}(H_{n_l}).
$$

 Using any constructive algorithm for computing the conformal radius \cite{BB}
we can approximate the $k$-th term of the sequence 
$$
R_k = r(B_{2^{-(k-1)}}(W_k),\al)+5\cdot 2^{4-\frac{l-1}{2}}
$$
By \lemref{rad-modulus1},  $R_k \ra r(c)$, and moreover, $\{R_k\}$ is a non-increasing sequence. Let 
$\rho_k$ be a dyadic approximation of $R_k$ that we compute so that $|\rho_k - R_k|<2^{-k}$. 
Let 
$$
r_k = \rho_k + 3 \cdot 2^{-k}. 
$$
Then $\{r_k\}$ is a computable sequence of dyadic numbers. We have
$$
\lim_{k\ra \infty} r_k = \lim _{k\ra \infty} \rho_k = \lim _{k\ra \infty} R_k = r(c),
$$
and for each $k$,
$$
r_k = \rho_k + 3 \cdot 2^{-k} \ge R_k + 2\cdot 2^{-k} \ge R_{k+1} + 4 \cdot 2^{-(k+1)} \ge
\rho_{k+1} + 3\cdot 2^{-(k+1)} = r_{k+1}.
$$
This shows that $r(c)$ is right-computable.
\end{proof}

We now want to prove the ``if" direction of Theorem \ref{thm:char}.
Given a computable sequence $\{r_n\}$ such that $r_n \searrow r$ we 
claim that we can construct a $c$ such that $r= r(c)$. It is easy 
to do in the case when $r=0$, as any parabolic parameter would do. From 
now on, we assume that $r>0$.  We will be using 
the following two lemmas. The first one is Lemma~3.1 of \cite{BY}, and the second
one is Lemma~4.2 of \cite{BBY2}.

\begin{lem}
\label{smlchg}
For any initial segment $I = [a_0, a_1, \ldots, a_n]$, write 
$\omega = [a_0, a_1, \ldots, a_n, 1, 1, 1, \dots]$. Then 
for any $\ve>0$, there is an $m>0$ and an integer $N$
such that if we write $\be = [a_0, a_1, \ldots, a_n, 1, 1, \ldots, 
1, N, 1, 1, \ldots]$, where the $N$ is located in the $n+m$-th 
position, then
$$
\Phi(\omega) + \ve < \Phi(\be) < \Phi(\omega) + 2 \ve. 
$$
\end{lem}

\begin{lem}
\label{notdeclem}
For  $\omega$ as above, for any $\ve>0$ there is an $m_0>0$, which 
can be computed from $(a_0, a_1, \ldots, a_n)$ and $\ve$,  
such that for any $m \ge m_0$, and for any tail 
$I = [a_{n+m}, a_{n+m+1}, \ldots]$ if we denote 
$$ \be^{I} = [a_1, a_2, \ldots, a_n, 1, 1, \ldots, 1, a_{n+m}, a_{n+m+1}, 
\ldots],$$ then
$$
\Phi(\be^{I}) > \Phi (\om) - \ve. 
$$
\end{lem}

Using Lemma \ref{bounded-type}, we can get a computable version of Lemmas \ref{smlchg} 
and \ref{notdeclem}. 

\begin{lem}
\label{smlchgr}
For any given initial segment $I = [a_0, a_1, \ldots, a_n]$ and $m_0>0$, write 
$\omega = [a_0, a_1, \ldots, a_n, 1, 1, 1, \dots]$. Then 
for any $\ve>0$, we can uniformly compute  $m>m_0$, an integer $t$ and an integer $N$
such that if we write $\be = [a_0, a_1, \ldots, a_n, 1, 1, \ldots, 
1, N, 1, 1, \ldots]$, where the $N$ is located in the $n+m$-th 
position, we have 
\beq
\label{rcond}
r(\omega) - 2 \ve < r(\be) < r(\omega) -  \ve,
\eeq
\beq
\label{phiinc}
\Phi(\be)> \Phi(\om),
\eeq
and for any
$$
\ga = [a_0, a_1, \ldots, a_n, 1, 1, \ldots, 1, N, 1, \ldots, 1, c_{n+m+t+1}, c_{n+m+t+2}, \ldots],
$$
\beq
\Phi(\ga)>\Phi(\om)-2^{-n}. 
\eeq
\end{lem}

\begin{proof}
We first show that such $m$ and $N$ {\em exist}, and then 
give an algorithm to compute them. By Lemma \ref{smlchg}
we can increase $\Phi(\omega)$ by any controlled amount 
by modifying one term arbitrarily far in the expansion. 

By Theorem \ref{phi-cont}, $f: \theta\mapsto \Phi(\theta)+\log r(\theta)$
extends to a continuous function. Hence for any $\ve_0$ there is a $\de$ such 
that $|f(x)-f(y)|<\ve_0$ whenever $|x-y|<\de$. In particular, 
there is an $m_1$ such that $|f(\be)-f(\om)| < \ve_0$ whenever $m \ge m_1$.

This means that if we choose $m$ large enough, a controlled increase 
of $\Phi$ closely corresponds to a controlled drop of $r$ by a corresponding 
amount, hence there are $m>m_0$ and $N$ such that 
\eref{rcond} holds. \eref{phiinc} is satisfied almost automatically.
The only problem is to {\em computably} find such $m$ and $N$. 

To this end, we apply Lemma \ref{bounded-type}. It implies that for any specific $m$ 
and $N$ we can compute $r(\be)$. This means that we can find the suitable 
$m$ and $N$, by enumerating all the pairs $(m,N)$ and exhaustively checking 
\eref{rcond} and \eref{phiinc} for all of them. We know that eventually we will find a 
pair for which \eref{rcond} and \eref{phiinc} hold. 

Finally, $t$ exists and can be computed by Lemma \ref{notdeclem}.
\end{proof}

We are now ready to give an algorithm for computing a rotation number $\theta$ 
for which $r(\theta)=\lim \searrow r_n$. $c$ is easily computed from $\theta$. 
The algorithm works as follows. On stage $k$ it produces a finite initial segment 
$I_k=[a_0, \ldots, a_{m_k}]$ such that the following properties are maintained:
\begin{enumerate}
\item 
$I_0$=[];
\item 
$I_k$ has at least $k$ terms, i.e. $m_k \ge k$;
\item 
for each $k$, $I_{k+1}$ is an extension of $I_k$;
\item 
for each $k$, denote $\ga_k = [I_k, 1, 1, \ldots]$, then $r_k+2^{-(k+1)}< \ga_k < r_k + 2^{-k}$; 
\item
for each $k$, $\Phi(\ga_k)> \Phi(\ga_{k-1})$;
\item 
for each $k$, for any extension 
$$
\be = [I_k, b_{m_k+1}, b_{m_k+2}, \ldots], 
$$
$\Phi(\be) > \Phi(\ga_k)-2^{-k}$. 
\end{enumerate}

The first three properties are very easy to assure. The last three are maintained 
using Lemma \ref{smlchgr}. By this Lemma we can decrease $r(\ga_{k-1})$ by any given 
amount (possibly in more than one step) by extending $I_{k-1}$ to $I_k$. Here we 
use the facts that the $r_k$'s are computable and non-increasing. 

Denote $\theta = \lim_{k\ra \infty} \ga_k$. The continued fraction expansion of 
$\theta$ is the limit of the initial segments $I_k$. 
This algorithm gives us at least one term of the continued fraction expansion 
of $\theta$ per iteration, hence we would need at most $O(n)$ iterations 
to compute $\theta$ with precision $2^{-n}$ (in fact, much fewer iterations 
would suffice).
It remains to prove that, in fact, $\theta$ is the rotation number we are looking 
for. 

\begin{lem}
The following equalities hold:
$$
\Phi(\theta) = \lim_{k \rightarrow \infty}\Phi(\gamma_k)~~~~\mbox{$~$and$~$}~~~~
r(\theta) = \lim_{k \rightarrow \infty}r(\gamma_k)=r.
$$
\end{lem}

\begin{proof}
By the construction, the limit $\theta = \lim \ga_k$ exists. We also 
know that the sequence $r(\ga_k)$ converges  to 
the number $r=\lim\searrow r_k$, and that the sequence $\Phi(\ga_k)$ is 
monotone non-decreasing, and hence converges to a value $\psi$ ({\em a priori}
we could have $\psi=\infty$). By the Carath{\'e}odory Kernel Theorem, we have $r(\theta) \ge r>0$, so 
$\Phi(\theta) < \infty$.
On the other hand, by the property we have maintained through 
the construction, we know that $\Phi(\theta)>\Phi(\ga_k)-2^{-k}$ for all $k$. 
Hence $\Phi(\theta) \ge \psi$.  In particular, $\psi < \infty$. 

From \cite{BC} we know that 
\beq
\label{BCCont}
\psi + \log r= \lim (\Phi(\ga_k) + \log r(\ga_k)) = \Phi(\theta) + \log r(\theta).
\eeq
Hence we must have $\psi = \Phi (\theta)$, and $r = r(\theta)$, which completes
the proof. 
\end{proof}

\section{Constructing non-computable Julia sets}
\label{sec:constructive}

Theorem \ref{thm:char} provides us with a tool for constructing 
explicit parameters $c$ for which $J_{z^2+c}$ is non-computable. 
In fact, for any right-computable $r$ which is not a computable 
number, we obtain such a set. 

In particular, we can give an explicit construction of a Julia set 
parameter, for which computing the Julia set with an arbitrarily
high precision would allow us to solve the {\em halting problem}. 

Let $R(x,t)$ be the predicate which evaluates to $1$ if and only 
if $x$ is a valid encoding of a Turing Machine which terminates 
after exactly $t$ steps. Then for any fixed $x$, $R(x,\bullet)$ evaluates
to $1$ at most once, and 
$$
H(x) = \exists t ~R(x,t)
$$
is the non-computable Halting predicate. Let $\{r_k\}$ be the following 
computable sequence of dyadic numbers
$$
r_k = \frac{1}{16} \left( 1 - \sum_{x=1}^{k}\sum_{t=1}^{k} 4^{-x}\cdot R(x,t) \right).
$$
$r_k$ is obviously a non-increasing sequence, and it converges to 
$$
r = \lim_{k \ra \infty} r_k = \frac{1}{16} \left( 1 - \sum_{x=1}^{\infty} 4^{-x}\cdot H(x) \right).
$$
Hence $r$ is a right-computable number, and there is a computable parameter
$c$ such that $r(c)=r$. Now note that computing $r=r(c)$ with precision 
$4^{-(x+3)}$ would require evaluating $H(1), H(2), \ldots, H(x)$, which is 
impossible.

\section{The complexity of $\theta$'s for which $J_\theta$ is non-computable}
\label{sec:complexity}

In this note we have demonstrated that one can explicitly compute a parameter 
$\theta$ such that the corresponding Julia set $J_\theta = J_{P_\theta}$ is 
non-computable. In fact, computing $J_\theta$ would allow us to solve the 
Halting Problem. It is a natural question to ask whether such a parameter 
$\theta$ can be computed fast. That is, how much time is required to produce
a $2^{-n}$-approximation of $\theta$. The algorithm for generating $\theta$
presented above does not yield any explicit bound on the computational complexity 
of $\theta$. In this section we show that under a reasonable conjecture
$\theta$ can be shown to be {\em poly-time} computable. That is, there is a 
polynomial $p(n)$ such that it takes at most $p(n)$ steps to produce a $2^{-n}$
approximation of $\theta$. The conjecture is a stronger version of Theorem \ref{phi-cont}.

\begin{conj}
\label{conj:mod}
The 1-periodic continuous function $f: \theta \mapsto \Phi(\theta) + \log r(\theta)$ has 
a computable modulus of continuity. In other words, there is a computable function 
$\mu : \NN \ra \NN$ such that $| f(\theta_1)-f(\theta_2)| < 2^{-n}$ whenever
$| \theta_1 - \theta_2|<2^{-\mu(n)}$. 
\end{conj}

Note that in Conjecture \ref{conj:mod} we do not put any restrictions on the
growth of $\mu$ besides it being computable. $\mu(n)$ can grow
extremely fast, as some computable integer functions do. 
In fact, the following much stronger conjecture has been put forward by Marmi,
Moussa, and Yoccoz in \cite{MMY}, and is generally supported by 
the available evidence:

\medskip
\noindent
{\bf Marmi-Moussa-Yoccoz Conjecture.} \cite{MMY} {\it The function $f:\theta\mapsto \Phi(\theta)+\log r(\theta)$ is H{\"o}lder of exponent $1/2$.}

\medskip
\noindent
Assuming the weaker Conjecture \ref{conj:mod}, we show that the $\theta$'s from 
Theorem \ref{thm:char} (and hence the $c$'s) can be made poly-time computable. 
Of course, the strengthening of Theorem \ref{thm:char} only comes in the ``if"
direction. We only state that direction here. 

\begin{thm}[{\bf Conditional}]
\label{thm:poly}
 Suppose Conjecture \ref{conj:mod} holds, and
suppose there is a computable sequence $r_1, r_2, \ldots$ of dyadic numbers such 
that 
\begin{itemize}
\item 
$\{r_i\}$ is non-increasing, $r_1 \ge r_2 \ge \ldots$, and 
\item 
$\lim_{i\ra \infty} r_i = r$. 
\end{itemize}
Then there is a {\em poly-time} computable $\theta$ (and hence a poly-time computable $c=c(\theta)$)
such that $r(c) = r$.
\end{thm}

\begin{proof}
The proof goes along the lines of the proof of the ``if" direction of Theorem \ref{thm:char}.
We outline the modifications made to the proof here and leave the details to the reader. 
The key difference is that in the proof of Theorem \ref{thm:char} we used Lemma \ref{smlchgr}
to perform a step in decreasing the conformal radius from $r(\ga_{k-1})$ to $r(\ga_k)$. 
The algorithm there is basically an exhaustive search, which of course could take much
more than a polynomial time in the precision of $\ga_k$ to compute. By assuming 
Conjecture \ref{conj:mod} we can deal with $\Phi(\ga_{k-1})$ and $\Phi(\ga_k)$ instead 
of the $r(\bullet)$'s. We have an explicit formula for $\Phi$ that converges well, and 
we can compute the continued fractions coefficient to make $\Phi(\ga_k)$ close to 
whatever we want relatively fast.

The step of going from $\ga_{k-1}$ to $\ga_k$ is as follows. First, we 
do the following computations:
\begin{itemize}
\item 
compute $d_k$ which is the ``drop" in $r$ we are trying to achieve; 
we want $d_k/2 < \log(r(\ga_{k-1})) - \log(r(\ga_k)) < d_k$;
\item 
compute using the function $\mu$ a value $\de_k$ such that $|f(x)-f(y)|<d_k/8$
whenever $|x-y|<\de_k$. 
\end{itemize}
We have no {\em a priori} bound on how long these computations would take, but 
we would still like to be computing $\theta$ in polynomial time. 
To achieve this, we use $1$'s in the continued fraction expansion of 
$\theta$ to ``pad" the computation. 

When asked about the value of $\theta$ with precision $2^{-n}$ which is higher 
than what the known terms of the expansion $[I_{k-1}]$ can provide, we do the following:
\begin{itemize}
\item 
try to compute $d_k$ and $\de_k$ as above, but run the computation for 
at most $n$ steps;
\item
if the computation does not terminate, output an answer consistent with 
the initial segment $[I_{k-1},\underbrace{1,1,\ldots,1}_{2 n}]$;
\item 
if the computation terminates in less than $n$ steps proceed as described 
below.
\end{itemize}
Note that so far the computation is polynomial in $n$. 
For some sufficiently large $n$ the computation will terminate in $n$ steps, 
at which point we will have computed $d_k$ and $\de_k$. We then add more $1$'s to the 
initial segment to assure that $|\ga_{k-1}-\ga_{k}|<\de_k$. 

Recall that our goal is to assure that 
$$d_k/2 < \log(r(\ga_{k-1})) - \log(r(\ga_k)) < d_k.$$
Withe the current initial segment for $\ga_k$ we have $|\ga_{k-1}-\ga_{k}|<\de_k$, and 
hence in the difference
$$
\log(r(\ga_{k-1})) - \log(r(\ga_k))  = 
\Phi(\ga_k) - \Phi(\ga_{k-1}) + (f(\ga_{k-1}) - f(\ga_k))
$$
the last term is bounded by $d_k/8$. 
This means that for the current step it suffices to increase $\Phi(\ga_k)$ relative to 
$\Phi(\ga_{k-1})$ by between $\frac{5}{8} d_k$ and $\frac{7}{8} d_k$. 

Let $M$ be the total length of $I_{k-1}$ and the $1$'s we have added, and let us extend the
continued fraction by putting $N\in\NN$ in the $M+1$-st term, and all $1$'s further.
Increasing $M$ if necessary, we can ensure an approximate equality 
$$
\Phi(\ga_k) \approx \Phi(\ga_{k-1}) + \al(N) \log N
$$
up to an error of $\frac{1}{32}d_k$.
Let $p_M/q_M$ be 
the $M$-th convergent of the resulting continued fraction. 
Recall that on an input $n$ we need to compute $\theta$ with 
precision $2^{-n}$ in time polynomial in $n$. If $2^{-n}> 1/\sqrt{q_M}$, 
then we do not need to know anything about $N$ to compute the required 
approximation. Suppose $2^{-n}<1/\sqrt{q_M}$, which means $n> \log q_M/2$. 
And we have time polynomial in $\log q_M$ to perform the computation. 

Note that $M = O(\log q_M)$. It is also not hard to see that $\al(N)<2^{M/2}$,
so in order to have a change by $\approx 3 d_k/4$ we must have $N>e^{\Om(2^{M/2})}$, 
hence by making $M$ sufficiently large (depending on the value of $d_k$), 
we can guarantee that $N>e^{2^{M/3}}$. This means that we can 
approximate $\al(N)$ with the truncated function $\Phi$ at the $M$-th
convergent of the continued fraction. Write $p_M/q_M = [a_1, a_2, \ldots, a_M]$, and denote
$$
\be = [a_1, a_2, \ldots, a_M] \cdot [a_2, a_3, \ldots, a_M] \cdot \ldots
\cdot [a_{M-1}, a_M] \cdot [a_M].
$$
Then $\be$ approximates $\al(N)$ within a very small relative error. In particular, we can 
assure that 
$$
\be \cdot \left( 1- \frac{1}{32}\right) < \al(N) < \be \cdot \left(1+\frac{1}{32}\right).
$$
In time polynomial in $\log q_m$ we can compute the exact 
expression for $\be$ using rational arithmetic: $\be = p/q$. 
Now we can estimate $N$ and write it as $e^{6 d_k/ 8\be}$ in time polynomial in $\log (q_M)$. 
From there we can continue by adding enough $1$'s to get $I_k$ and 
$\ga_k = [I_k,1,1,\ldots]$. By the construction, it would give us the necessary decrease
 in the value of $r(\ga_k)$. 
\end{proof}

 \end{document}